\documentclass[11pt, fullpage]{article}

\usepackage{amsfonts}
\usepackage{xr}
\usepackage{color}
\usepackage{graphicx,epsfig}
\usepackage{amsmath}
 \usepackage{color}
\usepackage[pdfencoding=unicode]{hyperref}
\usepackage[utf8]{inputenc}
\usepackage[T1]{fontenc}
\usepackage{svg}
\usepackage{pdfpages}

\newtheorem{theorem}{Theorem}[section]
\newtheorem{lemma}[theorem]{Lemma}

\newtheorem{definition}[theorem]{Definition}

\newtheorem{remark}[theorem]{Remark}

\numberwithin{equation}{section}

\def\R{\mathbb{R}}
\def\e{\varepsilon}

\def\tmm{{\theta_{m,\mu}}}

\def\pmm{{p_{m,\mu}}}

\begin{document}

\title{Optimization of the $L^{1}$ norm of the solution of a Fisher-KPP equation in the small diffusivity regime}

\date{}


\author{Gr\'egoire Nadin\footnote{Institut Denis Poisson, Universit\'e d'Orl\'eans, Universit\'e de Tours, CNRS, Orl\'eans, France ({\tt gregoire.nadin@cnrs.fr}).}}



\maketitle

\begin{abstract}
We investigate in the present paper the maximization problem for the functional $F_{\mu}(m):=\int_{0}^{1}(\tmm-cm)$, where $\tmm$ is the unique positive solution of $ -\mu \theta'' = \theta (m-\theta)$ in $(0,1)$, $ \theta'(0)=\theta'(1)=0$, and $0\leq m\leq \kappa$. We assume $c\in (1,3)$. It is already known that the $BV$ norms of maximizers of this functional blow up when the diffusivity $\mu$ tends to $0$. Here, we first show that the maximizers are always $BV(0,1)$. Next, we completely characterize the limit of the maximas of this functional as $\mu\to 0$, and we show that one can construct a quasi-maximizer which is periodic, in a sense, and with a $BV$ norm behaving like $1/\sqrt{\mu}$. Lastly, we prove that along a subsequence $\mu_{k}\to 0$, any maximizer of $F_{\mu_{k}}$ is periodic, in a sense.  
\end{abstract}

\noindent {\bf Key-words:} logistic diffusive equation, heterogeneous Fisher-KPP, regularity and symmetry in optimization problems. 

\smallskip

\noindent {\bf AMS classification.} 34B15, 34C14, 49J15, 49K15

\section{Introduction}

The aim of this article is to describe the maximizers of 
$$F_{\mu}(m):=\int_{0}^{1}(\tmm-cm)$$
where  $\theta=\tmm$ is the unique positive solution of 
\begin{equation} \label{eq:theta} -\mu \theta'' = \theta (m-\theta) \quad \hbox{ in } (0,1), \quad \theta'(0)=\theta'(1)=0,\end{equation}
which is well-defined for all $m\in \mathcal{A}\textcolor{black}{(0,1)}\backslash\{0\}$, where 
$$\mathcal{A}\textcolor{black}{(0,1)}:= \{ m\in L^{\infty}\textcolor{black}{(0,1)}, \ 0\leq m\leq \kappa \hbox{ a.e. } [0,1]\}.$$
We consider the distribution of resources $m$ in the set of admissible functions $\mathcal{A}\textcolor{black}{(0,1)}$. 

\subsection{Earlier works}

This problem has first been raised by Lou \cite{Lou2008} in a slightly different form, namely, he \textcolor{black}{investigated} the maximizers of 
$$G_{\mu}(m):=\int_{0}^{1}\tmm$$
over the $m\in\mathcal{A}\textcolor{black}{(0,1)}$, with $\int_{0}^{1}m=m_{0}$ prescribed. It is well-known that these two problems are dual in a sense. \textcolor{black}{Namely, if $\overline{m}$ maximizes $F_\mu$, then it is a maximizer of $G_\mu$ over $\mathcal{A}(0,1)$ with constraint $\int_0^1 m = \int_0^1 \overline{m}$. On the other hand, if $\overline{m}$ is a  maximizer of $G_\mu$ over $\mathcal{A}(0,1)$ with constraint $\int_0^1 m = m_0$, then there exists $c$, which depends on $m_0$, such that $\overline{m}$ maximizes $F_\mu$.} The first order optimality conditions were first derived in \cite{Ding2010}, where some numerics were also performed. 

The author addressed this problem with Mazari and Privat in \cite{Mazari2020}, and proved that only the two crenels $m=m_{cr,\ell}:=\kappa 1_{(0,\ell}$ and $m_{cr,\ell}(1-\cdot)$ are maximizers when $\mu$ is large enough. In particular, these maximizers are $BV(0,1)$ and {\em bang-bang}, that is, $m(x)\in \{0,\kappa\}$ for a.e. $x\in (0,1)$. Mazari \cite{Mazari} proved that the large diffusivity regime is related to the investigation of steady states for heterogeneous diffusive Lotka-Volterra competition systems. 

In parallel, Nagahara and Yanagida \cite{NagaharaYanagida} proved that the maximizers are bang-bang for all $\mu>0$, under the assumption that these maximizers are Riemann-measurable. 

The regime $\mu\to 0$ was investigated by Mazari and Ruiz-Balet \cite{MazariRuizBalet}, who proved that the $BV$ norms of the maximizers necessarily tend to $+\infty$ as $\mu\to 0$, meaning that the maximizers oscillate very fastly between $0$ and $\kappa$ when $\mu$ is small. They performed some precise numerics describing such a behaviour, including  multidimensional sets. 

This property was improved by Mazari, Privat and the author \cite{MNP}, who proved that the $BV$ norms of the maximizers blow-up at least as $C/\sqrt{\mu}$ as $\mu\to 0$ for some constant $C>0$.\textcolor{black}{We recall that the BV norm is defined as 
$$\|m\|_{BV(0,1)}:= \|m\|_{L^1(\Omega)}+\sup \big\{ \int_0^1 m(x)\phi'(x)dx, \ \phi\in \mathcal{C}^1_c(0,1), \ \|\phi\|_\infty\leq 1\big\}.$$
When $m$ is bang-bang, then it is equivalent to the perimeter of the set $\{m=\kappa\}$. 
}
They also proved in the same paper that the maximizers are always bang-bang, regardless of any regularity or large diffusivity assumption.  
The method used to derive this property is quite general and was used to derive bang-bang properties for a wide class of bilinear control problems by Mazari \cite{Mazaribang}. 

A discretized version of the problem, with discrete Laplacian, was investigated by Lou, Nagahara, and Yanagida in \cite{LouNagaharaYanagida}. In that case, they managed to fully describe the maximizers when $\mu\to 0$. These maximizers are close to a periodic function. However, the connection between this discrete problem and the continuous one is not clear in the small diffusivity regime $\mu\to 0$. 
Another discretized version of the equation, with equal diffusion rate between each patches, was investigated in \cite{LiangZhang}. 

Another related problem raised attention these last years: the maximization of the ratio $\int_{0}^{1}\tmm/\int_{0}^{1}m$ under the constraint $m\geq 0$, $m\not\equiv 0$ on the growth rate. Bai, He and Li \cite{BaiHeLi} proved that the supremum of this ratio is exactly $3$, and that a maximizing sequence $(m_{n})_{n}$ is the one concentrating to a Dirac mass at $x=1$. Inoue described the behaviour of $\tmm$ along such a sequence in \cite{Inoue}. This ratio is not bounded anymore in multidimensional domains \cite{InoueKuto}. 


\subsection{Statement of the results}

Let now come back to the maximization of 
$$F_{\mu}(m):=\int_{0}^{1}(\tmm-cm)\textcolor{black}{.}$$
We assume in the present paper that $c\in (1,3)$. 

If $c\geq 3$, then the unique maximizer of $F_{\mu}$ is $0$. Indeed, in that case, it has been proved in \cite{BaiHeLi} that $\int_{0}^{1}\tmm <3 \int_{0}^{1}m$ for any non-constant $m$. Hence, $F_{\mu}(m)<0$ as soon as $m$ is non constant for $c\geq 3$. As $\tmm\equiv m$ when $m$ is constant, one concludes that $0$ is the unique maximizer. 

If $c\leq 0$, then clearly $m\equiv \kappa$ is the unique global maximizer. 

Thus, only the case $c\in (0,1]$ remains relevant and is not covered by the present paper. We explain in Remark \ref{rmq:c1} the main obstacles in trying to extend the present method to $c\in (0, 1]$. 

\bigskip

We start with a regularity result on the maximizers. 

\begin{theorem} \label{theorem:BV} Assume that $c\in (1,3)$. Let $\overline{m}_{\mu}$ \textcolor{black}{be} a maximizer of $F_{\mu}$ and assume that $\overline{m}_{\mu}\not\equiv 0$ and $\overline{m}_{\mu}\not\equiv \kappa$. Then the function $\theta_{\overline{m}_{\mu},\mu}'$ admits a finite number of zeros, and $\overline{m}_{\mu}$ admits exactly one jump from $0$ to $\kappa$ \textcolor{black}{or from $\kappa$ to $0$} between each of these zeros. In particular, $\overline{m}_{\mu}$ is $BV(0,1)$. 
\end{theorem}

We will prove in Lemma \ref{lem:sup>0}, using a result of \cite{BaiHeLi}, that $\max_{\mathcal{A}\textcolor{black}{(0,1)}}F_{\mu}>0$ when $\mu$ is small enough. Hence, the hypotheses $\overline{m}_{\mu}\not\equiv 0$ and $\overline{m}_{\mu}\not\equiv \kappa$ are satisfied for $\mu$ small enough. 

\bigskip

We now introduce a notion of functions of particular interest. Such functions appeared in the numerics performed in \cite{MazariRuizBalet}. 

\begin{definition}\label{def:ksym}
We say that a function $m\in L^{\infty}(0,1)$ is {\em $k-$symmetric}, for $k\in \mathbb{N}\backslash \{0\}$, if there exists a function $m_{0}\in L^{\infty}(0,1)$, such that for all $l\leq \left \lfloor{k/2}\right \rfloor $:
$$m(x):= \left\{ \begin{array}{lcll}
m_{0} (\textcolor{black}{k}y)&\hbox{ if }& x=\displaystyle\frac{2l}{k}+y,& \quad y\in [0,\frac{1}{k}),\\
m_{0}(1-\textcolor{black}{k}y) &\hbox{ if }& x=\displaystyle\frac{2l+1}{k}+y,& \quad y\in [0,\frac{1}{k}).\\
\end{array}\right. $$

\end{definition}

An example of a $k-$symmetric function is given at the bottom of Figure \ref{fig:mmu} below. If $k=2l$ is even, then $m$ is periodic, in the sense that it could be written $m(x)=m_{per}(lx)$ for some $1-$periodic even function $m_{per}$. Hence, the notion of $k-$symmetry is somehow an extended notion of periodicity. 

Let us denote for all $\ell \in [0,1]$ the crenel distribution as:
$$m_{cr,\ell}(x):= \left\{ \begin{array}{rcl} \kappa &\hbox{ if }& x\in [0,\ell],\\
0 &\hbox{ if }& x\in (\ell,1].\\ \end{array}\right. $$
Define 
$$G(\mu):=\max_{\ell\in [0,1]}\textcolor{black}{F_{\mu}}(m_{cr,\ell}) \quad \hbox{ and } \quad \textcolor{black}{\overline{\mu}_{inf}} := \inf\{\mu>0, G(\mu)=\sup_{\mu'>0}G(\mu')\}$$

\begin{theorem} \label{theorem:mu0} Assume that $c\in (1,3)$. Then $\textcolor{black}{\overline{\mu}_{inf}}$ is well-defined and positive. For all $\mu>0$, there exists $\hat{m}_{\mu}\in \mathcal{A}(0,1)$ such that 
 $\|\hat{m}_{\mu}\|_{BV(0,1)}=k_{\mu}$ as $\mu\to 0$, $\hat{m}_{\mu}$ is $k_{\mu}$-symmetric with pattern $m_{cr,\ell}$ for some $\ell \in (0,1)$, and 
$$\max_{m'\in \mathcal{A}(0,1)}F_{\mu}(m')-2\kappa \sqrt{\mu/\textcolor{black}{\overline{\mu}_{inf}}}\leq F_{\mu}(\hat{m}_{\mu})\leq \max_{m'\in \mathcal{A}(0,1)}F_{\mu}(m').$$
Lastly, $\max_{m'\in \mathcal{A}(0,1)}F_{\mu}(m')\to G(\textcolor{black}{\overline{\mu}_{inf}})=\textcolor{black}{\sup_{\mu>0}G(\mu)}$ as $\mu \to 0$. 
\end{theorem}

\begin{remark}
We leave as an open problem the conjecture that $G$ admits a unique maximizer. This would imply the convergence of $\mu \textcolor{black}{k_{\mu}^{2}}$ as $\mu\to 0$ to this unique maximizer. Liang and Lou \cite{LiangLou} provided an example of growth rate $m$ for which $\mu \mapsto F_{\mu}(m)$ admits at least two maximas. However, the growth rates considered in the present paper are quite different from the one considered in \cite{LiangLou}, which was the perturbation of a constant function, and we thus still believe that \textcolor{black}{the maximizer $ \textcolor{black}{\overline{\mu}_{inf}}$ is unique} here. 
\end{remark}

We do not know if the maximizers are always $k-$symmetric for $\mu$ small. In the numerics performed in \cite{MazariRuizBalet}, the maximizers did not always look like $k-$symmetric functions. However, when the diffusivity is well scaled, in a sense, we are able to prove such a $k-$symmetry. 

\begin{theorem}\label{theorem:equality}  Assume that $c\in (1,3)$. 
If one can write $\mu=\overline{\mu}/k^{2}$ for some $k\in \mathbb{N}\backslash \{0\}$ and $\overline{\mu}$ a maximizer of $\mu>0\mapsto G(\mu)$, then any maximizer $\overline{m}_{\mu}$ of $\textcolor{black}{F_{\mu}}$ is $K-$symmetric with pattern $m_{cr,\ell}$ for some $\ell>0$ and $K\in \mathbb{N}\backslash \{0\}$. Moreover, if $G$ admits a unique maximizer $\overline{\mu}$, then $K=k$. 
\end{theorem}


\section{$BV$ regularity of the maximizers}

The aim of this section is to prove Theorem \ref{theorem:BV}. 

Consider the Hamiltonian \textcolor{black}{
defined for all $x\in [0,1], \theta, \theta', m \in \R, \eta \in \R^{2}$ by
$$H(x,\theta, \theta', \eta, m):=  \eta_{1}\theta'-\frac{1}{\mu}\eta_{2}\theta (m-\theta)-\frac{1}{\mu}(\theta -cm)$$
and the cost 
$$C(m) := \int_{0}^{1}(-\tmm+cm)= -F_{\mu}(m).$$}

\textcolor{black}{
The Hamiltonian is related to $\tmm$ through the equation:
$$\frac{d}{dt}\begin{pmatrix} \tmm \\ \theta^{'}_{m,\mu}\\ \end{pmatrix} = 
\partial_{\eta}H(x,\tmm, \theta^{'}_{m,\mu}, \eta, m) $$
and there exists a solution $\eta$ of the adjoint equation
$$\frac{d}{dt}\begin{pmatrix} \eta_{1} \\ \eta_{2}\\ \end{pmatrix} =- \begin{pmatrix}\partial_{\theta}H(x,\tmm, \theta^{'}_{m,\mu}, \eta, m)\\
\partial_{\theta'}H(x,\tmm, \theta^{'}_{m,\mu}, \eta, m) \\ \end{pmatrix}= \begin{pmatrix} \frac{1}{\mu}\big((m-2\tmm)\eta_{2}+1\big)\\
-\eta_{1} \\ \end{pmatrix} $$}

\textcolor{black}{Defining $\pmm:= \eta_{2}$, we find that it is indeed} the unique solution of 
\begin{equation} \label{eq:p}
-\mu p'' - (m-2\tmm)p=1 \quad \hbox{ in } (0,1), \quad p'(0)=p'(1)=0. \end{equation}

The Pontryagin maximum principle (see for instance \textcolor{black}{Theorem 3 of Chapter 5 of \cite{LeeMarkus}}) yields that a minimizer $\overline{m}_{\mu}$ of the cost (that is, a maximizer of $F_{\mu}$) \textcolor{black}{minimizes 
$$m\in [0,\kappa]\mapsto H\big(x,\tmm (x), \theta^{'}_{m,\mu}(x), -p^{'}_{m,\mu}(x),\pmm(x), m)$$
for all $x\in [0,1]$, that is:}
\begin{equation} \label{eq:bathtub} \overline{m}_{\mu} (x)= \left\{ \begin{array}{lll}
\kappa &\hbox{ if }& c<\pmm(x)\tmm(x),\\
0 &\hbox{ if }& c>\pmm(x)\tmm(x).\\
\end{array} \right.\end{equation}

We will now denote $m=\overline{m}_{\mu}$, $\theta=\tmm$ and $p=\pmm$ when there is no ambiguity in order to enlight the notations, and
$$\Phi:= \theta p.$$

\textcolor{black}{Lastly, t}he Hamiltonian $H$ is constant along the trajectories:
\begin{equation} \label{eq:PMP}
\mu p'\theta'+p\theta (m-\theta)+\theta - cm = cste \hbox{ in } (0,1).\end{equation}

%
%
%
%

The next lemmas provide as a by-product an alternative proof of the bang-bang property of maximizers, which could  also be obtained by duality from Theorem I of \cite{MNP}.

\begin{lemma} \label{lem:zeros} Assume that $c>1$. 
Let $0\leq a<1$ \textcolor{black}{be} such that $\theta'(a)=0$ and $\Phi (a)<c$. Then 
$$b:= \min \{x\in (a,1], \ \theta'(x)=0\}$$
is well-defined, $\theta'>0$ in $(a,b)$, $\Phi$ only crosses $c$ once in $(a,b]$, and $\Phi (b)>c$. 

Similarly, if $\Phi (a)>c$, then $b$ is still well-defined and $\theta'<0$ in $(a,b]$, $\Phi$ only crosses $c$ once in $(a,b)$, and $\Phi (b)<c$. 
\end{lemma}

{\bf Proof.} 

We just prove the first part, the other one being proved similarly. As $\Phi$ is continuous, we know that $\Phi<c$ in a right neighborhood $(a,a+\e)$ of $a$ and thus $m=0$ in $(a,a+\e)$. It follows that $\mu \theta''=\theta^{2}> 0$ in $(a,a+\e)$ (since $\theta>0$ due to $m\not\equiv 0$) and thus, as $\theta'(a)=0$, $\theta'>0$ in $(a,a+\e)$. As $\textcolor{black}{\theta'(1)}=0$, $b$ is well-defined and strictly larger than $a$. 

Next, assume $\Phi (b)\leq c$ by contradiction. We know that if $\Phi< c$ in $(a,b)$, then one would have $m=0$ a.e. and $\theta''>0$ a.e. in $(a,b)$, contradicting $\theta'(b)=0$. Hence, there exists an interval $[x_{0},y_{0}]\subset [a,b]$ such that $\Phi \geq c$ in this interval. Moreover, we could assume that $\Phi (y_{0})=p(y_{0})\theta (y_{0})=c$ and as $\Phi (a)=p(a)\theta (a)<c$ , one gets from (\ref{eq:PMP}):
$$\begin{array}{rcl} 
(1-c)\theta(a) &< & -p(a)\theta^{2}(a)+\theta (a)\\
&&\\
&=& \mu p'(y_{0})\theta'(y_{0})+c (m(y_{0})-\theta(y_{0}))+\theta(y_{0}) - cm(y_{0})\\
&&\\
&=& \mu p'(y_{0})\theta'(y_{0})+(1-c)\theta(y_{0}). \\
\end{array}$$

%

As $\theta'>0$ over $(a,b)$, one has $\theta (y_{0})>\theta (a)$ and thus, as $c>1$, one gets 
$p'(y_{0})>0$. This implies $\Phi'(y_{0})>0$, a contradiction since $\Phi (y_{0})=c$ and $\Phi \geq c$ in $(x_{0},y_{0})$. We have also proved that $\Phi$ only crosses $c$ once in $(a,b]$, otherwise, there exist $(x_{0},y_{0})$ as above and we could conclude similarly. 
\hfill $\Box$

\begin{lemma} \label{lem:zerosdeg} Assume that $c>1$. 
Then for all $a\in [0,1]$ such that $\theta'(a)=0$, one has $\Phi (a)\neq c$. 
\end{lemma}

{\bf Proof.} 

Assume that there exists  $a\in [0,1]$ such that $\theta'(a)=0$, and $\Phi (a)= c$. \textcolor{black}{Let us prove by contradiction that this implies that 
for all $\tilde{a}\in [0,1]$ such that $\theta'(\tilde{a})=0$, one has $\Phi (\tilde{a})= c$. Hence, consider $\tilde{a}\in [0,1]$ such that $\theta'(\tilde{a})=0$ and $\Phi (\tilde{a})\neq c$.}
If $\tilde{a}<a$, as any such $\tilde{a}$ is isolated by Lemma \ref{lem:zeros}, one can assume that $\tilde{a}$ is the largest one satisfying this property. But then, either $\Phi(\tilde{a})>c$ and then $b:= \min \{x\in (a,L], \ \theta'(x)=0\}$ satisfies $\Phi (b)<c$, which contradicts the definition of $\tilde{a}$ since $b$ is either $a$ or another point $x$ such that $\theta'(x)=0$ and $\Phi (x)=c$. Similarly, if $\Phi(\tilde{a})<c$, then $\Phi (b)>c$ and we also reach a contradiction. If $\tilde{a}>a$, we reach a similar contradiction by assuming that $\tilde{a}$ is the smallest one satisfying this property and using again Lemma \ref{lem:zeros}. Hence, we have proved by contradiction that for all $\tilde{a}$ such that $\theta'(\tilde{a})=0$, one has $\Phi (\tilde{a})=c$. In particular $\Phi (0)=\Phi (1)=c$. 

Next, assume that there exists an interval $(x_{0},y_{0})$ such that $\Phi>c$ on this interval. We could assume that 
$\Phi (x_{0})=\Phi (y_{0})=c$. The function $\theta'$ does not vanish on $(x_{0},y_{0})$ since $\Phi>c$. 
Assume first that $\theta'>0$ on $(x_{0},y_{0})$. We define $a=\max \{x\in [0,x_{0}], \theta'(x)=0\}$ and 
$b=\min \{x\in [0,x_{0}], \theta'(x)=0\}$. Then $\Phi (a)=\Phi(b)=c$ since $\theta'(a)=\theta'(b)=0$ and the same computations as in the proof of Lemma \ref{lem:zeros} yield
$$(1-c)\theta (a)=\mu p'(y_{0})\theta'(y_{0})+(1-c)\theta (y_{0}).$$
It follows that $p'(y_{0})>0$ since $\theta(y_{0})>\theta(a)$ and $c>1$, a contradiction since $\Phi (y_{0})=c$ and $\Phi>c$ in $(x_{0},y_{0})$. If there exists an interval $(x_{0},y_{0})$ such that $\Phi<c$ on this interval, we reach a contradiction similarly. 

Hence, we have proved by contradiction that if there exists  $a\in [0,1]$ such that $\theta'(a)=0$, and $\Phi (a)= c$, then $\Phi\equiv c$ on $(0,1)$. This is a contradiction with Lemma \ref{lem:zeros}. 

\hfill $\Box$

{\bf Proof of Theorem \ref{theorem:BV}.} 

By Lemmas \ref{lem:zeros} and \ref{lem:zerosdeg}, we know that all the zeros of $\theta'$ in $[0,1]$ are isolated. Hence, as $[0,1]$ is compact, $\theta'$ only admits a finite number of zeros. Hence, by Lemma \ref{lem:zeros}, $\Phi$ only crosses $c$ a finite number of times. By characterization (\ref{eq:bathtub}), $m$ is $BV(0,1)$ and admits exactly one jumps between each zeros of $\theta'$. 

\hfill $\Box$

\begin{remark} \label{rmq:c1} We do not know if Theorem \ref{theorem:BV} still holds when $0<c\leq 1$. The main difference when  $c<1$ is that the constant function $m(x)=\kappa$ satisfies the first order optimality conditions. Indeed, in that case $\theta\equiv \kappa$, $p\equiv 1/\kappa$ and $\Phi\equiv 1\geq c$ almost everywhere. Hence, Lemmas \ref{lem:zeros} and \ref{lem:zerosdeg} cannot hold for any functions $(\theta,p,m)$ satisfying $(\ref{eq:theta})$, $(\ref{eq:p})$, $(\ref{eq:bathtub})$ and $(\ref{eq:PMP})$. However, $m\equiv \kappa$ might only be a local maximizer of $F_{\mu}$, not a global one, and maybe Lemmas \ref{lem:zeros} and \ref{lem:zerosdeg} still hold for global maximizers. We leave this possible extension as an open problem. 
\end{remark}

\section{Existence of a quasi-maximizer}

The aim of this section is to prove Theorem \ref{theorem:mu0}. 

In all this section,  we will specify the dependence of $F_{\mu}$ with respect to the interval considered. That is, we define 
$$F_{\mu}^{(a,b)}(m):=\int_{a}^{b}(\theta-cm)$$
where $0\leq a<b\leq 1$ and $\theta$ is the unique positive solution of 
\begin{equation} \label{eq:ab} -\mu \theta'' = \theta (m-\theta) \quad \hbox{ in } (a,b), \quad \theta'(a)=\theta'(b)=0.\end{equation}

\subsection{Construction of the quasi-maximizer $\hat{m}_{\mu}$}\label{sec:construction}

We consider $m=\overline{m}_{\mu}$ such that $F^{(0,1)}_{\mu}(m)=\max_{m'\in \mathcal{A}(0,1)}F_{\mu}^{(0,1)}(m')$.
By Theorem \ref{theorem:BV}, we denote by $(a_{i})_{1\leq i\leq N+1}$ the zeros of $\theta'$, with $a_{0}=0$ and $a_{N+1}=1$, and we know that $m$ only jumps from $0$ to $\kappa$ \textcolor{black}{or from $\kappa$ to $0$} once in each interval $(a_{i},a_{i+1})$. In other words, 
$$\textcolor{black}{\|m_{|(a_{i},a_{i+1})}\|_{BV(a_{i},a_{i+1})}}=1.$$

\begin{lemma} \label{lem:decomp}
One has for all $i\geq 1$:
$$F^{(0,1)}_{\mu}(m)=\sum_{i= 0}^{N}F^{(a_{i},a_{i+1})}_{\mu}(m).$$
\end{lemma}

{\bf Proof.}

As $\theta'(a_{i})=\theta'(a_{i+1})=0$, the solution $\theta$ of (\ref{eq:ab}), with $a=a_{i}$ and $b=a_{i+1}$ is just $\theta$ restricted to $(a_{i},a_{i+1})$ by uniqueness. Hence, $F^{(a_{i},a_{i+1})}_{\mu}(m)=\int_{a_{i}}^{a_{i+1}}(\theta-cm)$. The decomposition $F^{(0,1)}_{\mu}(m)=\sum_{\textcolor{black}{i= 0}}^{N}F^{(a_{i},a_{i+1})}_{\mu}(m)$ follows. 

\hfill $\Box$

\bigskip

We now define 
$$A_{i}:= \frac{1}{a_{i+1}-a_{i}}F^{(a_{i},a_{i+1})}_{\mu}(m).$$

It follows from Lemma \ref{lem:decomp} that 
$$F^{(0,1)}_{\mu}(m)=\sum_{i= 0}^{N}(a_{i+1}-a_{i})A_{i}\leq \max_{0\leq i\leq N}A_{i}.$$
Let $i_{\mu}$ \textcolor{black}{be} such that $\max_{0\leq i \leq N}A_{i}= A_{i_{\mu}}$. One obtains 
$$ F^{(0,1)}_{\mu}(m)\leq A_{i_{\mu}}\leq  \frac{1}{a_{i_{\mu}+1}-a_{i_{\mu}}}\max_{m'\in \mathcal{A}(a_{i_{\mu}},a_{i_{\mu}+1})}F_{\mu}^{(a_{i_{\mu}},a_{i_{\mu}+1})}(m').$$

\begin{figure}[h]
    \centering
    \includegraphics[width=10cm]{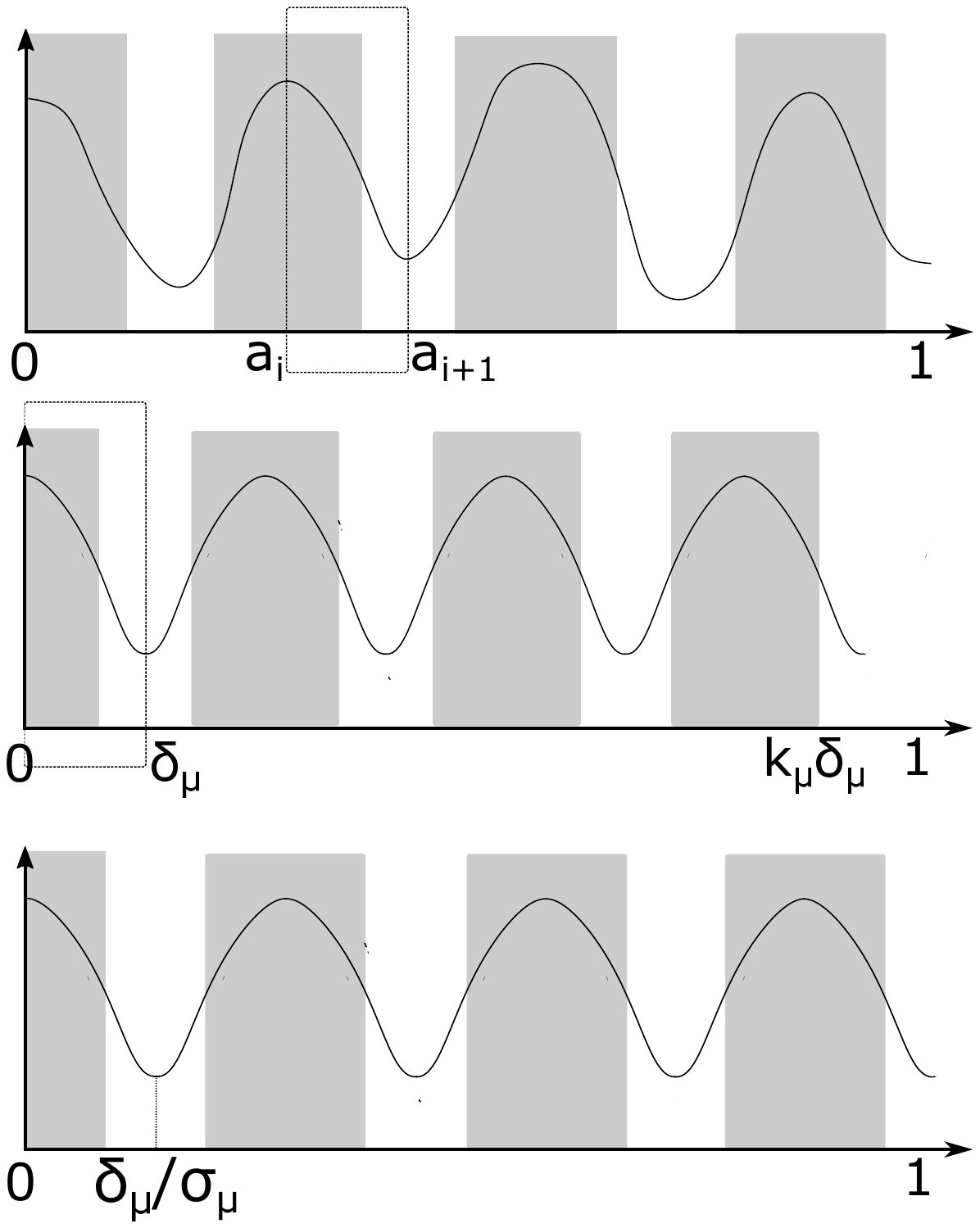}
\caption{Construction of $\hat{m}_{\mu}$. First step: select the interval $(a_{i},a_{i+1})$ 
    maximizing $A_{i}$. Second step: repeat periodically this pattern. Third step: stretch it so that it ends at $x=1$.}
\label{fig:mmu}
\end{figure}

Let 
$$1=k_{\mu}(a_{i_{\mu}+1}-a_{i_{\mu}})+r_{\mu}, \quad \hbox{ with } \quad 0\leq r_{\mu}<a_{i_{\mu}+1}-a_{i_{\mu}}.$$
We now construct a $k_{\mu}-$symmetric function  
from $m_{|(a_{i_{\mu}},a_{i_{\mu}+1})}$
by symmetrization and dilation. Namely, we define $\delta_{\mu}:= a_{i_{\mu}+1}-a_{i_{\mu}}$ and 

$$\hat{m}_{\mu}(x):= \left\{ \begin{array}{lcll}
m\big(\sigma_{\mu} y+a_{i_{\mu}}\big) &\hbox{ if }& x=2l\frac{\delta_{\mu}}{\sigma_{\mu}}+y,& \quad y\in [0,\frac{\delta_{\mu}}{\sigma_{\mu}}),\\
m\big(a_{i_{\mu}+1}-\sigma_{\mu} y\big) &\hbox{ if }& x=(2l+1)\frac{\delta_{\mu}}{\sigma_{\mu}}+y,& \quad y\in [0,\frac{\delta_{\mu}}{\sigma_{\mu}}),\\
\end{array}\right. $$
for all $l$ such that $l\leq  \left \lfloor{k_{\mu}/2}\right \rfloor$, and where $\sigma_{\mu}:=k_{\mu}\delta_{\mu}=1-r_{\mu}$. This construction is described in Figure \ref{fig:mmu}.

\begin{lemma} \label{lem:dilatation}
One has for all $\lambda>0$, $\tilde{m}\in \mathcal{A}(0,1)$:
$$F_{\mu}^{(0,1)}(\tilde{m})=\lambda F_{\frac{\mu}{\lambda^{2}}}^{(0,1/\lambda)}\big(\tilde{m}(\lambda \cdot)\big).$$
\end{lemma}

{\bf Proof.}

Let $\theta:= \theta_{\tilde{m},\mu}$ and $\theta_{\lambda}(x):= \theta (\lambda x)$.  One has 
$$-\frac{\mu}{\lambda^{2}}\Delta \theta_{\lambda} (x)= \theta (\lambda x)\big( \tilde{m}(\lambda x)-\theta (\lambda x)\big)= \theta_{\lambda}( x)\big( \tilde{m}(\lambda x)-\theta_{\lambda}( x)\big) \hbox{ on } (0,1/\lambda).$$
Hence, 
$$F_{\frac{\mu}{\lambda^{2}}}^{(0,1/\lambda)}\big(\tilde{m}(\lambda \cdot)\big) = \int_{0}^{1/\lambda}(\theta_{\lambda}-c\tilde{m}(\lambda \cdot))= \frac{1}{\lambda}\int_{0}^{1}(\theta-c\tilde{m})=\frac{1}{\lambda}F^{(0,1)}_{\mu}(\tilde{m}).$$

\hfill $\Box$

\begin{lemma} \label{lem:ksymm}
Assume that $\tilde{m}\in \mathcal{A}(0,1)$ is a $k-$symmetric function with pattern $\tilde{m}_{0}$. Then for all $\mu>0$,
$$F_{\mu}^{(0,1)}(\tilde{m})=F_{k^{2}\mu}(\tilde{m}_{0}).$$
\end{lemma}

{\bf Proof.}

Clearly, $F_{\mu}^{(0,1)}(\tilde{m}) = \sum_{i=0}^{k-1}F_{\mu}^{(\frac{i}{k}, \frac{i+1}{k})}(\tilde{m})$. Moreover, by \textcolor{black}{ definition of $\tilde{m}$, one has $F_{\mu}^{(\frac{i}{k}, \frac{i+1}{k})}(\tilde{m})=F_{\mu}^{(0, \frac{1}{k})}\big(\tilde{m}_{0}(k\cdot)\big)$ if $i$ is even, $F_{\mu}^{(\frac{i}{k}, \frac{i+1}{k})}(\tilde{m})=F_{\mu}^{(0, \frac{1}{k})}\big(\tilde{m}_{0}(1-k\cdot)\big)$ if $i$ is uneven. As $F_{\mu}^{(0, \frac{1}{k})}\big(\tilde{m}_{0}(1-k\cdot)\big)=F_{\mu}^{(0, \frac{1}{k})}\big(\tilde{m}_{0}(k\cdot)\big)$ by symmetry, one has: } 
$$F_{\mu}^{(0,1)}(\tilde{m}) = k F_{\mu}^{(0, \frac{1}{k})}(\tilde{m})= F_{k^{2}\mu}^{(0,1)}\big(\tilde{m}(\cdot/k)\big) =F_{k^{2}\mu}^{(0,1)}(\tilde{m}_{0})$$
by Lemma \ref{lem:dilatation}. 
\hfill $\Box$

\begin{lemma}\label{lem:calculhatm}
One has 
$$F_{\mu}^{(0,1)}(\hat{m}_{\mu})=\frac{1}{a_{i_{\mu}+1}-a_{i_{\mu}}}F_{\mu(1-r_{\mu})^{2}}^{(a_{i_{\mu}},a_{i_{\mu}+1})}(m).$$
\end{lemma}

{\bf Proof.}

Let $\tilde{m}(x):= \hat{m}_{\mu}(x/\sigma_{\mu})$. This is the function corresponding to the second step described in Figure \ref{fig:mmu}. One has by Lemma \ref{lem:dilatation}:
$$\sigma_{\mu} F_{\mu}^{(0,1)}\big(\hat{m}_{\mu}\big)=F_{\mu \sigma_{\mu}^{2}}^{(0,\sigma_{\mu})}(\tilde{m})=\sum_{j=0}^{k_{\mu}-1}F_{\mu \sigma_{\mu}^{2}}^{(j\delta_{\mu},(j+1)\delta_{\mu})}(\tilde{m}).$$
Next, it is clear that for all $j$, $F_{\mu \sigma_{\mu}^{2}}^{(j\delta_{\mu},(j+1)\delta_{\mu})}(\tilde{m})=F_{\mu \sigma_{\mu}^{2}}^{(a_{i_{\mu}},a_{i_{\mu}+1})}(m)$ by symmetry and definition of $\hat{m}_{\mu}$. 
Hence, 
$$\sigma_{\mu} F_{\mu}^{(0,1)}\big(\hat{m}_{\mu}\big)=k_{\mu}F_{\mu \sigma_{\mu}^{2}}^{(a_{i_{\mu}},a_{i_{\mu}+1})}(m)$$
and the conclusion follows from $\sigma_{\mu}=k_{\mu}\delta_{\mu}=1-r_{\mu}$.
\hfill $\Box$

%
%
%

\begin{lemma}\label{lem:estF'}
For all $\tilde{m}\in \mathcal{A}(0,1)$, $\mu\mapsto F_\mu (\tilde{m})$ is of class $\mathcal{C}^{1}$ and 
$$\partial_{\mu }F_{\mu }(\tilde{m}) \leq \kappa/\mu.$$
\end{lemma}

{\bf Proof.} 

The regularity follows from classical arguments and one has  $\partial_{\mu }F_{\mu }(\tilde{m})=\int_{0}^{1}\dot{\theta}$, where $\dot{\theta}$ is the unique solution of 
$$-\mu \dot{\theta}''-(\tilde{m}-2\theta_{\tilde{m},\mu})\dot{\theta} = \theta_{\tilde{m},\mu}'' \hbox{ in } (0,1), \quad \dot{\theta}'(0)=\dot{\theta}'(1)=0.$$
Straightforward computations yield that $\theta_{\tilde{m},\mu}/\mu$ is a supersolution of this equation and thus the weak maximum principle yields 
$$\dot{\theta}\leq \theta_{\tilde{m},\mu}/\mu\leq \kappa /\mu,$$
from which the conclusion follows. 
 \hfill $\Box$

\begin{lemma}
For all $\tilde{m}\in \mathcal{A}(0,1)$, $\mu>0$ and $r\in (0,1)$, one has 
$$F_{\mu(1-r)^{2}}^{(0,1)}(\tilde{m})-F_{\mu}^{(0,1)}(\tilde{m})\geq -2\kappa r.$$
\end{lemma}

{\bf Proof.}

There exists $\xi \in \big( \mu (1-r)^{2},\mu\big)$ such that 
$$F_{\mu(1-r)^{2}}^{(0,1)}(\tilde{m})-F_{\mu}^{(0,1)}(\tilde{m})=\partial_{\mu}F_{\xi}^{(0,1)}(\tilde{m})\big(\mu(1-r)^{2}-\mu\big)\geq \textcolor{black}{-2\mu\partial_{\mu}F_{\xi}^{(0,1)}r.}$$
 \textcolor{black}{Hence, by Lemma \ref{lem:estF'}, one gets 
$$F_{\mu(1-r)^{2}}^{(0,1)}(\tilde{m})-F_{\mu}^{(0,1)}(\tilde{m})\geq  -2\kappa r.$$}
\hfill $\Box$

\bigskip

Gathering all the previous estimates, we have thus obtained the following  intermediate result. 

\begin{lemma}\label{prop:reste} There exists $\hat{m}_{\mu}\in \mathcal{A}(0,1)$ such that 
 $\|\hat{m}_{\mu}\|_{BV(0,1)}=k_{\mu}$,  $\hat{m}_{\mu}$ is $k_{\mu}-$symmetric with pattern $m_{cr,\ell_{\mu}}$ for some $\ell_{\mu}\in [0,1]$ and 
$$\max_{m'\in \mathcal{A}(0,1)}F_{\mu}(m)-2\kappa \delta_{\mu}\leq F_{\mu}(\hat{m}_{\mu})\leq \max_{m'\in \mathcal{A}(0,1)}F_{\mu}(m).$$
\end{lemma}

\subsection{Estimates on $\mu/\delta_{\mu}^{2}$}

The aim of this section is to prove Proposition \ref{prop:accu}, which, combined with Proposition \ref{prop:reste}, ends the proof of Theorem \ref{theorem:mu0}.

\begin{lemma} \label{prop:accu} Define $k_{\mu}$, $\delta_{\mu}$ and $\ell_{\mu}$ as in the previous section. Consider a sequence $(\mu_{k})_{k}$ such that $\lim_{k\to+\infty}\mu_{k}=0$, 
$(\mu_{k}/\delta_{\mu_{k}}^{2})_{k}$ converges to $\mu_{0} \geq 0$, and $(\ell_{\mu_{k}})_{k}$ converges to $\ell_{0}$ as $k\to +\infty$. Then 
$$F_{\mu_{0}}^{(0,1)}(m_{cr,\ell_{0}})\geq F_{\mu}^{(0,1)}(m_{cr,\ell}) \quad \hbox{ for all } \ell \in [0,1], \mu>0.$$
In particular, 
$$\textcolor{black}{\overline{\mu}_{inf}}\leq \liminf_{\mu\to 0}\mu/\delta_{\mu}^{2}\leq \limsup_{\mu\to 0}\mu/\delta_{\mu}^{2}\leq \textcolor{black}{\overline{\mu}_{sup}}$$ 
where, if we let $G(\mu):=\max_{\ell\in [0,1]}F_{\mu}^{(0,1)}(m_{cr,\ell})$, one defines 
$$\textcolor{black}{\overline{\mu}_{inf}} = \min \{\mu>0, G(\mu)=\sup_{\mu'>0}G(\mu')\} \hbox{ and } \textcolor{black}{\overline{\mu}_{sup}} = \max \{\mu>0, G(\mu)=\sup_{\mu'>0}G(\mu')\},$$
and these two quantities are positive and finite. 

As a consequence, $\max_{m'\in \mathcal{A}(0,1)}F_{\mu}(m')\to G(\textcolor{black}{\overline{\mu}_{inf}})=\sup G$ as $\mu \to 0$. 
\end{lemma}

\begin{figure}[!h]
    \centering
    \includegraphics[width=12cm]{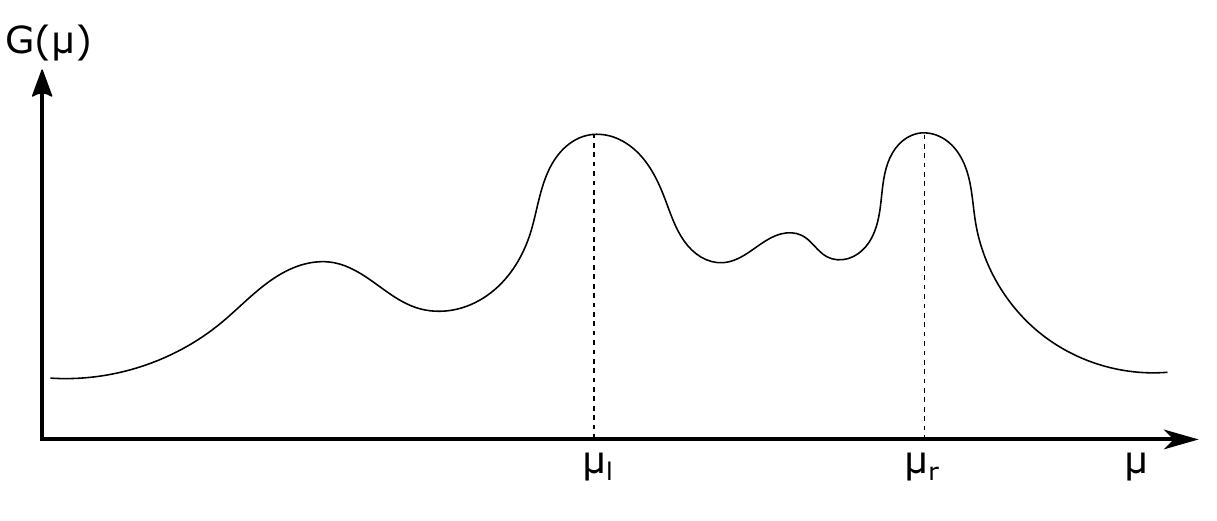}
    \caption{An illustration of the definitions of $\textcolor{black}{\overline{\mu}_{inf}}$ and $\textcolor{black}{\overline{\mu}_{sup}}$.}
    \label{fig:Gmu}
\end{figure}


Let first compare the supremum $F^{(0,1)}_{\mu}(m)$ with an appropriate function.
Consider $\ell\in [0,1]$, $k\in\mathbb{N}\backslash \{0\}$  and $\tilde{m}(x)$ a $k-$symmetric function with pattern $m_{cr,\ell}$. Gathering all the previous inequalities and using Lemma \ref{lem:ksymm}, we have:
\begin{equation} \label{ineq:kl}A_{i_{\mu}}=F^{(0,1)}_{\mu/\delta_{\mu^{2}}}(m_{cr,\ell_{\mu}})\geq F^{(0,1)}_{\mu}(m)\geq F^{(0,1)}_{\mu}(\tilde{m})=F^{(0,1)}_{\mu k^{2}}(m_{cr,\ell}).\end{equation}

We now need to make sure that $\mu/\delta_{\mu^{2}}$ is not too large nor too small.

\begin{lemma}\label{lem:mugrand}
There exists a constant $C>0$ such that for all $m\in \mathcal{A}\textcolor{black}{(0,1)}$ and $\mu>0$, one has 
$$\big\|\theta_{m,\mu}-\int_{0}^{1}m\big\|_{L^{\infty}(0,1)}\leq C\kappa^{3/2}/\sqrt{\mu}.$$
\end{lemma}

{\bf Proof.} 

Multiplying the equation satisfied by $\tmm$ by $\tmm$ and integrating by parts, one gets 
$$\mu \int_{0}^{1}(\tmm')^{2}= \int_{0}^{1}\tmm^{2}(m-\tmm)\leq \kappa^{3}.$$
It follows from the $(\infty,2)$ Poincar\'e inequality \cite{Mazya} that there exists a constant $C>0$ such that 
$$\big\|\tmm - \int_{0}^{1}\tmm\big\|_{L^{\infty}(0,1)}\leq \frac{C\kappa^{3/2}}{6\sqrt{\mu}}.$$
On the other hand, integrating the equation satisfied by $\tmm$, one gets
$$\int_{0}^{1}\tmm m = \int_{0}^{1}\tmm^{2}$$
and thus, writing $\tmm = \int_{0}^{1}\tmm + \eta$:
$$\int_{0}^{1}\tmm \big( \int_{0}^{1}(\tmm - m)\big) = \int_{0}^{1} \eta m - 2 \int_{0}^{1}\tmm \int_{0}^{1}\eta - \int_{0}^{1}\eta^{2}\leq \int_{0}^{1} \eta m - 2 \int_{0}^{1}\tmm \int_{0}^{1}\eta.$$
As $\int_{0}^{1}\tmm \geq \int_{0}^{1}m$, one gets 
$$0\leq \int_{0}^{1}(\tmm - m)\leq \|\eta\|_{L^{\infty}(0,1)}  +2 \int_{0}^{1}|\eta|\leq C\kappa^{3/2}/(2\sqrt{\mu}),$$
from which the conclusion follows. 
 \hfill $\Box$

\begin{lemma}\label{lem:sup>0}
For all $c<3$, one has 
$$\limsup_{\mu\to 0}\sup_{\ell \in [0,1]}F_{\mu}^{(0,1)} (m_{cr,\ell})=\sup_{\mu>0}\sup_{\ell \in [0,1]}F_{\mu}^{(0,1)} (m_{cr,\ell})>0.$$
\end{lemma}

{\bf Proof.}

The identity between the $\limsup$ and the $\sup$ could be proved exactly as in \cite{MazariRuizBalet}. 

Next, we consider the sequence $m_{cr,\e}$, that is
$$m_{cr,\e}(x):= \left\{ \begin{array}{rcl} \kappa &\hbox{ if }& x\in [0,\e],\\
0 &\hbox{ if }& x\in (\e,1].\\ \end{array}\right. $$
It has been proved in \cite{BaiHeLi} that 
$$F_{\sqrt{\e}}^{(0,1)}(\frac{1}{\kappa \e}m_{cr,\e}) \to 3 \int_{0}^{1}\frac{m_{cr,\e}}{\kappa\e}\quad \hbox{ as } \e \to 0.$$
Hence, we could consider $c'\in (c,3)$ and $\e$ small enough such that 
$$F_{\sqrt{\e}}^{(0,1)}(\frac{1}{\kappa \e}m_{cr,\e})+c\geq 3 - c'+c,$$
where we have used $\int_{0}^{1}m_{cr,\e}=\kappa\e$. 
Next, one easily checks by considering $\theta_{B}:= B\tmm$ that for all $\mu>0$, $B>0$ and $m\in L^{\infty}(0,1)$, $m\geq 0$, one has 
$$F_{B\mu}^{(0,1)}(Bm)=B F_{\mu}^{(0,1)}(m).$$
Hence, 
$$F_{\kappa\e^{3/2}}^{(0,1)}(m_{cr,\e}) \geq (3 - c')\e\kappa>0,$$
which ends the proof. 
\hfill $\Box$

\begin{lemma}\label{lem:mu>0}
There exists $M>0$ such that $1/M\leq \mu/\delta_{\mu}^{2}\leq M$ for all $\mu>0$. 
\end{lemma}

{\bf Proof.} 

Assume first that there exists a sequence $(\mu_{j})_{j}$ such that $\lim_{j\to \infty} \mu_{j}/\delta_{\mu_{j}}^{2}=0$. We can assume, up to extraction, that there exists $\ell_{0}\in [0,1]$ such that $\lim_{j\to +\infty}\ell_{\mu_{j}}=\ell_{0}$. It has been proved in \cite{MazariRuizBalet} that for all $M>0$, $\|\tmm-m\|_{L^{1}(0,1)}\to 0$ as $\mu\to 0$ uniformly on function $m\in\mathcal{A}(0,1)$ such that $\|m\|_{BV(0,1)}\leq M$. 
In particular, as crenels have $BV-$norms equal to $1$, one has $\|\theta_{m_{cr,\ell},\mu}-m_{cr,\ell}\|_{L^{1}(0,1)}\to 0$ as $\mu\to 0$ uniformly with respect to $\ell \in [0,1]$. 

We thus obtain 
$$A_{i_{\mu_{j}}}=F^{(0,1)}_{\mu_{j}/\delta_{\mu_{j}^{2}}}(m_{cr,\ell_{\mu_{j}}})\to (1-c)\int_{0}^{1}m_{cr,\ell_0}= \kappa (1-c)\ell_{0}\leq 0 \hbox{ as } j\to +\infty.$$
As $$\lim_{j\to +\infty}A_{i_{\mu_{j}}}\geq \limsup_{j\to +\infty}F_{\mu_{j}}^{(0,1)}(\overline{m}_{\mu_{j}})= \limsup_{\mu\to 0}F_{\mu}^{(0,1)} (\overline{m}_{\mu})=\sup_{\mu>0}F_{\mu}^{(0,1)} (\overline{m}_\mu)>0$$
by Lemma \ref{lem:sup>0}, one reaches a contradiction. 

Next, if  there exists a sequence $(\mu_{j})_{j}$ such that $\lim_{j\to \infty} \mu_{j}/\delta_{\mu_{j}}^{2}=+\infty$, then Lemma \ref{lem:mugrand} yields 
$$A_{i_{\mu_{j}}}=F^{(0,1)}_{\mu_{j}/\delta_{\mu_{j}^{2}}}(m_{cr,\ell_{\mu_{j}}})\simeq (1-c)\int_{0}^{1}m_{cr,\ell_{\mu_{j}}}=\kappa (1-c)\ell_{\mu_{j}}\leq 0 \hbox{ as } j\to +\infty, $$
leading to a contradiction again. 
 \hfill $\Box$

\bigskip

{\bf Proof of Proposition \ref{prop:accu}.}

Let $(\mu_{j})_{j}$ \textcolor{black}{be} such that $\mu_{j}\to 0$,  $\mu_{j}/\delta_{\mu_{j}}^{2}\to \mu_{0}\in [1/M,M]$ and $\ell_{\mu_{j}}\to \ell_{0}$ as $j\to +\infty$. 

Let $\eta>0$, $\ell>0$ and write $\sqrt{\eta}= k_{j}\sqrt{\mu_{j}}+r_{j}$, for some $k_{j}\in \mathbb{N}\backslash \{0\}$ and $0\leq r_{j}<\sqrt{\mu_{j}}$, so that $\lim_{j\to +\infty}\mu_{j}k_{j}^{2}=\eta$. Hence, it follows from (\ref{ineq:kl}) that 
$$F^{(0,1)}_{\mu_{0}}(m_{cr,\ell_{0}})\geq \limsup_{j\to +\infty}\sup_{\mathcal{A}\textcolor{black}{(0,1)}}F^{(0,1)}_{\mu_{j}}\geq \liminf_{j\to +\infty}\sup_{\mathcal{A}\textcolor{black}{(0,1)}}F^{(0,1)}_{\mu_{j}}\geq F^{(0,1)}_{\eta}(m_{cr,\ell}).$$
Hence, as this is true for any $\ell>0$ and $\eta>0$, all these inequalities are indeed equalities, which implies that 
$\mu_{0}$ is a maximizer of $G$ and that 
$$\lim_{\mu\to 0}\sup_{\mathcal{A}\textcolor{black}{(0,1)}}F^{(0,1)}_{\mu}= \sup_{\eta} G(\eta).$$
%

Lastly, it follows from the same arguments as in the proof of Lemma \ref{lem:mu>0} that $\textcolor{black}{\overline{\mu}_{inf}}>0$ and $\textcolor{black}{\overline{\mu}_{sup}}<\infty$. 

\hfill $\Box$

\section{The equality case}

The aim of this section is to prove Theorem \ref{theorem:equality}. We begin with a characterization of $k-$symmetric functions. 
We denote again $m=\overline{m}_{\mu}$, $\theta=\tmm$ and $p=\pmm$ when there is no ambiguity.

\begin{lemma}\label{lem:ksym}
Assume that there exists $a\in (0,1)$ such that $\theta'(a)=0$ and $p'(a)=0$. Then there exists $k\in \mathbb{N}$, $k\geq 1$, such that $m$ is $k-$symmetric. 
\end{lemma}

{\bf Proof.} 

By using the change of variable $x':=1-x$ if necessary, we can always assume that $a\in (0,1/2]$. 
Consider the symmetrized function $m_{s}(x):= m(2a-x)$ if $x\in [a,2a]$, $m(x)$ if $x\in [0,a]$, and define similarly $\theta_{s}$, and $p_{s}$. These functions satisfy $(\ref{eq:theta})$ and $(\ref{eq:p})$ on $(0,2a)$. 

We know that $\Phi (a)\neq c$ by Lemma \ref{lem:zerosdeg}.
We can assume that $\Phi (a)>c$, the other case being treated similarly. Let 
$$x_{0}:=\inf \{x\in (0,a), \Phi >c \hbox{ in } (x,a) \} \hbox{ and }  y_{0}:=\inf \{x\in (a,2a), \Phi (x)>c\}.$$
By continuity of $\Phi$, $x_{0}<a<y_{0}$. 
One has $m(x)=\kappa$ for all $x\in (a,y_{0})$ and $m_{s}(x)=\kappa$ for all $x\in (a,2a-x_{0})$. 
Let $z_{0}:=\min (y_{0},2a-x_{0})>a$. Then $\theta$ and $\theta_{s}$ both satisfy
$$\mu \theta''+\theta (\kappa - \theta)=0 \hbox{ in } (a,z_{0}).$$
Moreover, $\theta (a)=\theta_{s}(a)$ and $\theta'(a)=\theta_{s}'(a)(=0)$ by definition. Hence, $\theta\equiv \theta_{s}$ in $(a,z_{0})$. Similarly,  $p\equiv p_{s}$ in $(a,z_{0})$. As $\Phi=\theta p$, it follows that $x_{0}=2a-y_{0}$. 

Next, by Lemma \ref{lem:zeros}, we can define $b:= \min \{x\in (a,1], \ \theta'(x)=0\}$, $\theta'<0$ in $(a,b)$, and $\Phi$ only crosses $c$ once in $(a,b]$. Hence, if $z_{0}<b$, then $\Phi<c$ and $m=0$ in $(z_{0},b)$. As $\theta'(z_{0})=\theta_{s}'(z_{0})(=0)$, $\theta (z_{0})=\theta_{s}(z_{0})$, and $\mu\theta''=\theta^{2}$ in $(z_{0},b)$, one gets $\theta\equiv \theta_{s}$ in $(z_{0},b)$. Similarly,  $p\equiv p_{s}$ in $(z_{0},b)$. Moreover, by Lemma \ref{lem:zeros}, either $b=1$ or $\Phi (b)<c$. In particular, $\theta'(b)=0$ and $p'(b)=0$. We can thus iterate until $b=1$ or $b=2a$, thus proving that $\theta\equiv \theta_{s}$ and $p=p_{s}$ on $(0,b)$. Hence, $m=m_{s}$ on $(0,b)$

Considering now $\tilde{m}_{s}$ the symmetrized function with respect to $x=b$, we can prove using the same method that $m=\tilde{m}_{s}$ on $(0,2b)$ if $2b\leq 1$, on $(0,1)$ otherwise. Going on iterating, we conclude that $m$ is $k-$symmetric.
\hfill $\Box$

\bigskip

\noindent {\bf Proof of Theorem \ref{theorem:equality}.}

We define $(A_{i})_{i}$ and $i_{\mu}$ as in Section \ref{sec:construction}, and we have already proved that 
$$ F^{(0,1)}_{\mu}(\overline{m}_{\mu})\leq A_{i_{\mu}}=\frac{1}{a_{i_{\mu}+1}-a_{i_{\mu}}}F_{\mu}^{(a_{i_{\mu}},a_{i_{\mu}+1})}(\overline{m}_{\mu}).$$
Moreover, as $\overline{m}_{\mu}$ only jumps once between $0$ and $\kappa$, we could assume (up to the change of variable $x'=1-x$) that 
$$\overline{m}_{\mu}=\left\{ \begin{array}{lcl} \kappa &\hbox{ if }& x\in (a_{i_{\mu}},a_{i_{\mu}}+\delta_{\mu} \ell),\\
 0 &\hbox{ if }& x\in (a_{i_{\mu}}+\delta_{\mu} \ell, a_{i_{\mu}+1}),\\
 \end{array}\right. $$
  for some $\ell \in [0,1]$.  
 It follows from Lemma \ref{lem:dilatation} that 
 $$A_{i_{\mu}}= F_{\mu/\delta_{\mu}^{2}}^{(0,1)}(m_{cr,\ell}).$$

Take $\overline{\mu}$ as in the hypothesis of Theorem \ref{theorem:equality}, that is, $\mu =\overline{\mu}/k^{2}$ and $\overline{\mu}$ maximizes $G$. Take $\overline{\ell}$ such 
$F_{\overline{\mu}}^{(0,1)}(m_{cr,\overline{\ell}})=G(\overline{\mu})$. 
Define $\tilde{m}_{\mu}$ a $k-$symmetric function with pattern $m_{cr,\overline{\ell}}$. 

Then, as $\overline{m}_{\mu}$ is a maximizer in $\mathcal{A}(0,1)$, one has 
$$F_{\mu}^{(0,1)}(\overline{m}_{\mu})\geq F_{\mu}^{(0,1)}(\tilde{m}_{\mu}).$$
But Lemma \ref{lem:ksymm} yields:
$$F_{\mu}^{(0,1)}(\tilde{m}_{\mu})= F_{\mu k^{2}}^{(0,1)}(m_{cr,\overline{\ell}})=G(\overline{\mu}).$$

Gathering all these inequalities, we have proved that 
$$F_{\mu/\delta_{\mu}^{2}}^{(0,1)}(m_{cr,\ell})\geq F_{\mu}^{(0,1)}(\overline{m}_{\mu})\geq F_{\mu}^{(0,1)}(\tilde{m}_{\mu})\geq F_{\mu k^{2}}^{(0,1)}(m_{cr,\overline{\ell}})=G(\overline{\mu}).$$
As $F_{\mu/\delta_{\mu}^{2}}^{(0,1)}(m_{cr,\ell})\leq G(\mu/\delta_{\mu}^{2})\leq G(\overline{\mu})$ by definitions of $G$ and $\overline{\mu}$, this chain of inequalities is indeed an equality. Moreover, 
this chain of equalities also implies $A_{i}=A_{i_{\mu}}$ for all $i$.

In particular, $A_{0}=F_{\mu}^{(0,1)}(\overline{m}_{\mu})$, which means that $m|_{(0,a_{1})}$ maximizes $F_{\mu}^{(0,a_{1})}$. Define $q$ the adjoint function on $(a_{0},a_{1})$, where we remind to the reader that $a_{0}:=0$, that is, $q$ is the solution of 
\begin{equation} \label{eq:q}-\mu q''-(m-2\theta)q = 1 \hbox{ on } (0,a_{1}), \quad q'(0)=q'(a_{1})=0.\end{equation}
As $m|_{(0,a_{1})}$ maximizes $F_{\mu}^{(0,a_{1})}$, one has $m(x)=\kappa$ if $q(x)\theta(x)>c$ and $m(x)=0$ if $q(x)\theta(x)<c$. We know from Theorem \ref{theorem:BV} that $m$ only jumps once from $0$ to $\kappa$ in $(0,a_{1})$. Let $\ell$ \textcolor{black}{be} the point where the value of $m$ changes. Then $q(\ell)\theta(\ell)=c$. but we also know that $p(\ell)\theta(\ell)=c$. Hence, $q(\ell)=p(\ell)$. Moreover, $q'(0)=p'(0)=0$ and these two functions both satisfy (\ref{eq:q}). Let $z:=p-q$. One has 
$$-\mu z''-(m-2\theta)z=0 \hbox{ in } (0,\ell), \quad z'(0)=0, \ z(\ell)=0.$$
Moreover, if $y:=z/\theta$, then 
$$-\mu y''-2\mu \frac{\theta'}{\theta}y'+\theta y=0 \hbox{ in } (0,\ell), \quad y'(0)=0, \ y(\ell)=0.$$
It follows from the elliptic maximum principle that $y\equiv 0$ and thus $z\equiv 0$ and $p\equiv q$ in $(0,\ell)$. This identity extends to $(0,a_{1})$ by the Cauchy-Lipschitz theorem, and in particular, $p'(a_{1})=0$. It follows from Lemma \ref{lem:ksym} that $m$ is K-symmetric for some $K$.

It follows from Lemma \ref{lem:ksymm} that:
$$F_{\mu}^{(0,1)}(m)=F_{\mu K^{2}}^{(0,1)}(m_{cr,\ell})=G(\overline{\mu}).$$
Hence, if $G$ admits a unique maximizer $\overline{\mu}$, then $\mu K^{2}=\overline{\mu}=\mu k^{2}$ by hypothesis, and thus $k=K$. \hfill $\Box$

\begin{remark}
If $G$ admits several maximizers, it follows from the above arguments that $\mu K^{2}$ and $\mu k^{2}$ are both maximizers of $G$, \textcolor{black}{possibly} with $k\neq K$ two positive integers. This seems very unlikely and we believe that $G$ admits a unique maximizer. 
\end{remark}

%
%
%

\section{Discussion and open problems}

The main restriction of the present paper is $c>1$ and we leave as an open problem the case $c\in (0,1]$. We have already discussed about this hypothesis above. The reader could notice that, even without extending Theorem \ref{theorem:BV}, many results would extend if one could prove that the BV norm of the maximizers is uniformly bounded with respect to $\mu$ between two critical points of $\tmm$. Also, it would be good to reformulate the hypothesis $c>1$ in terms of an hypothesis on $m_{0}$ when considering the maximization problem for $G_{\mu}$ on the functions $m\in\mathcal{A}$ such that $\int_{0}^{1}m=m_{0}$. 

About the $BV$ norm of $\overline{m}_{\mu}$, we have proved in \cite{MNP} that it is bounded from below by $C/\sqrt{\mu}$. In the present paper, we have proved that it is bounded and that a function with BV norm $k_{\mu}$ of order $1/\delta_{\mu}$ is a quasi-maximizer. We leave as an open problem to show that $\textcolor{black}{\overline{\mu}_{inf}}=\textcolor{black}{\overline{\mu}_{sup}}$ defined in Proposition \ref{prop:accu} are equal, which would imply that $k_{\mu}\simeq \sqrt{\textcolor{black}{\overline{\mu}_{inf}}/\mu}$. Also, we do not know if one can prove a bound from above of order $1/\sqrt{\mu}$ on the $BV$ norm of original maximizer $\overline{m}_{\mu}$, or the convergence of $\|\overline{m}_{\mu}\|_{BV(0,1)}\sqrt{\mu}$ as $\mu\to 0$. 

Next, we have constructed a quasi-maximizer $\hat{m}_{\mu}$ using $\overline{m}_{\mu}$, but we were not able to show that these two functions are close in a sense. 

Lastly, the present method only works in dimension $1$, and the multidimensional framework remains open. The numerics displayed in \cite{MazariRuizBalet} indicate that the optimizers in multidimensional domains might be particularly irregular when $\mu\to 0$, despite some patterns seem to emerge.

\subsection*{Acknowledgments} The author would like to thank Idriss Mazari for fruitful discussions on this topic and relevant comments on its preliminary version. 

\bibliographystyle{amsplain}
\bibliography{BiblioFrag}

\end{document}